\documentclass[letterpaper, 10 pt, conference]{ieeeconf}  

\IEEEoverridecommandlockouts                              
\usepackage{flushend}

\usepackage{amsmath} 
\usepackage{amssymb}  
\usepackage{adjustbox}
\usepackage[font={footnotesize ,it}]{caption}
\usepackage{dsfont}
\usepackage{cite}
\usepackage{color}%

\newtheorem{thm}{\bf{Theorem}}[section]
\newtheorem{remark}{\bf{Remark}}[section]

\newcommand{\nn}{\nonumber}
\newcommand{\vepsilon}{\boldsymbol{\epsilon}}

\title{\LARGE \bf Binary Signaling under Subjective Priors and Costs as a Game}

\author{Serkan~Sar{\i}ta\c{s}, Sinan~Gezici and Serdar~Y\"uksel
\thanks{This research was supported in part by the Natural Sciences and Engineering Research Council (NSERC) of Canada.}
\thanks{S. Sar{\i}ta\c{s} and S. Gezici are with the Department of Electrical and Electronics Engineering, Bilkent University, 06800, Ankara, Turkey.
        {\tt\small \{serkan,gezici\}@ee.bilkent.edu.tr}}%
\thanks{S. Y\"uksel is with the Department of Mathematics and Statistics, Queen's University, Kingston, Ontario, Canada, K7L 3N6.
        {\tt\small yuksel@mast.queensu.ca}}%
}

\begin{document}

\maketitle
\thispagestyle{empty}
\pagestyle{empty}

\begin{abstract}
Many decentralized and networked control problems involve decision makers which have either misaligned criteria or subjective priors. In the context of such a setup, in this paper we consider binary signaling problems in which the decision makers (the transmitter and the receiver) have subjective priors and/or misaligned objective functions. Depending on the commitment nature of the transmitter to his policies, we formulate the binary signaling problem as a Bayesian game under either Nash or Stackelberg equilibrium concepts and establish equilibrium solutions and their properties. In addition, the effects of subjective priors and costs on Nash and Stackelberg equilibria are analyzed. It is shown that there can be informative or non-informative equilibria in the binary signaling game under the Stackelberg assumption, but there always exists an equilibrium. However, apart from the informative and non-informative equilibria cases, under certain conditions, there does not exist a Nash equilibrium when the receiver is restricted to use deterministic policies. For the corresponding team setup, however, an equilibrium typically always exists and is always informative. Furthermore, we investigate the effects of small perturbations in priors and costs on equilibrium values around the team setup (with identical costs and priors), and show that the Stackelberg equilibrium behavior is not robust to small perturbations whereas the Nash equilibrium is.
\end{abstract}

\section{INTRODUCTION}

In many decentralized and networked control problems, decision makers have either misaligned criteria or have subjective priors, which necessitates solution concepts from game theory. For example, detecting attacks, anomalies, and malicious behavior with regard to security in networked control systems can be analyzed under a game theoretic perspective, see e.g., \cite{sandberg2015cyberphysical, teixeira2015secure, networkSecurity,cyberSecuritySmartGrid, dan2010stealth,varshneyHypothesisGame,detectionGameAdversary, gupta2010optimal, gupta2012dynamic, basar1985complete}. 

In this paper, we consider signaling games that refer to a class of two-player games of incomplete information in which an informed decision maker (transmitter or encoder) transmits information to another decision maker (receiver or decoder) in the hypothesis testing context. In the following, we first provide the preliminaries and introduce the problems considered in the paper, and present the related literature briefly.

\subsection{Preliminaries}

Consider a binary hypothesis-testing problem:
\begin{align*}
\mathcal{H}_0 : Y = S_0 + N \;,\nn\\
\mathcal{H}_1 : Y = S_1 + N \;,
\end{align*}
where $Y$ is the observation (measurement) that belongs to observation set $\Gamma=\mathbb{R}$, $S_0$ and $S_1$ denote the deterministic signals under hypothesis $\mathcal{H}_0$ and hypothesis $\mathcal{H}_1$, respectively, and $N$ represents a Gaussian noise; i.e., $N \sim \mathcal{N} (0,\sigma^2)$. In the Bayesian setup, it is assumed that the prior probabilities of $\mathcal{H}_0$ and $\mathcal{H}_1$ are available, which are denoted by $\pi_0$ and $\pi_1$, respectively, with $\pi_0+\pi_1=1$.

In the conventional Bayesian framework, the aim of the
receiver is to design the optimal decision rule (detector) based
on $Y$ in order to minimize the Bayes risk, which is defined as \cite{Poor}
\begin{align}
r(\delta) = \pi_0 R_0(\delta) + \pi_1 R_1(\delta) \;,
\label{eq:riskEq}
\end{align}
where $\delta(\cdot)$ is the decision rule, and $R_i(\cdot)$ is the conditional risk of the decision rule when hypothesis $\mathcal{H}_i$ is true for $i\in\{0,1\}$. In general, a decision rule corresponds to a partition of the observation set $\Gamma$ into two subsets $\Gamma_0$ and $\Gamma_1$, and the decision becomes $\mathcal{H}_i$ if the observation $y$ belongs to $\Gamma_i$, where $i\in\{0,1\}$.

The conditional risks in \eqref{eq:riskEq} can be calculated as
\begin{align}
R_i(\delta) = C_{0i}\mathsf{P}_{0i} + C_{1i} \mathsf{P}_{1i} \;,
\label{eq:condRiskEq}
\end{align}
for $i\in\{0,1\}$, where $C_{ji}\geq0$ is the cost of deciding for
$\mathcal{H}_j$ when $\mathcal{H}_i$ is true, and $\mathsf{P}_{ji}=\mathsf{P}(Y\in\Gamma_j|\mathcal{H}_i)$ represents the conditional probability of deciding for $\mathcal{H}_j$ given that $\mathcal{H}_i$ is true, where $i,j\in\{0,1\}$ \cite{Poor}.

It is well-known that the optimal decision rule $\delta$ that minimizes the Bayes risk is the following likelihood ratio test (LRT):
\begin{align}
\delta : \Bigg\{ \pi_1 (C_{01}-C_{11}) p_1(y) \overset{\mathcal{H}_1}{\underset{\mathcal{H}_0}{\gtreqless}} \pi_0 (C_{10}-C_{00})p_0(y) \;,
\label{eq:lrt}
\end{align}
where $p_i(y)$ represents the probability density function (PDF)
of $Y$ under $\mathcal{H}_i$, where $i\in\{0,1\}$ \cite{Poor}.

If the transmitter and the receiver have the same objective function specified by \eqref{eq:riskEq} and \eqref{eq:condRiskEq}, then the signals can be designed to minimize the Bayes risk corresponding to the decision rule in \eqref{eq:lrt}. This leads to a conventional formulation which has been studied intensely in the literature \cite{Poor,KayDetection}. On the other hand, in order to reflect the different perspectives of the players, the transmitter and the receiver can have non-aligned Bayes risks. In particular, let $C^t_{ji}$ and $C^r_{ji}$ represent the cost values from the perspective of the transmitter and the receiver, respectively, where $i,j\in\{0,1\}$. Also let $\pi_i^t$ and $\pi_i^r$ for $i\in\{0,1\}$ denote the priors from the perspective of the transmitter and the receiver, respectively, with $\pi_0^j+\pi_1^j=1$, where $j\in\{t,r\}$. Here, from transmitter's and receiver's perspectives, the priors must be mutually absolutely continuous with respect to each other; i.e., $\pi_i^t\pi_i^r=0\Leftrightarrow\pi_i^t=\pi_i^r=0$ for $i\in\{0,1\}$. This condition assures that the impossibility of any hypothesis holds for both the transmitter and the receiver simultaneously. The aim of the transmitter is to perform the optimal design of signals $\mathbf{S}=\{S_0,S_1\}$ to minimize his Bayes risk; whereas, the aim of the receiver is to determine the optimal decision rule $\delta$ over all possible decision rules $\Delta$ to minimize his Bayes risk.

The Bayes risks are defined as follows for the transmitter and the receiver:
\begin{align*}
r^j(\mathbf{S},\delta) = \pi_0^j R^j_0(\mathbf{S},\delta) + \pi_1^j R^j_1(\mathbf{S},\delta)\;,
\end{align*}
where
\begin{align*}
R^j_i(\mathbf{S},\delta) = C^j_{0i} \mathsf{P}_{0i} + C^j_{1i} \mathsf{P}_{1i}\;,
\end{align*}
for $i\in\{0,1\}$ and $j\in\{t,r\}$. Here, the transmitter performs the optimal signal design problem under the power constraint below:
\begin{align*}
\mathbb{S}\triangleq\{\mathbf{S}:\lVert S_0 \rVert ^2 \leq P_0 \,,\; \lVert S_1 \rVert ^2 \leq P_1\} \;,
\end{align*}
where $P_0$ and $P_1$ denote the power limits.

In the simultaneous-move game, the encoder and the decoder announce their policies at the same time, and a pair of policies $(\mathbf{S}^*, \delta^*)$ is said to be a {\bf Nash equilibrium} \cite{basols99} if
\begin{align}
\begin{split}
r^t(\mathbf{S}^*, \delta^*) &\leq r^t(\mathbf{S}, \delta^*) \quad \forall \,\mathbf{S} \in \mathbb{S}\;, \\
r^r(\mathbf{S}^*, \delta^*) &\leq r^r(\mathbf{S}^*, \delta) \quad \forall \,\delta \in \Delta\;.
\label{eq:nashEquilibrium}
\end{split}
\end{align}
As noted from the definition in \eqref{eq:nashEquilibrium}, under the Nash equilibrium, each individual player chooses an optimal strategy given the strategies chosen by the other players.

However, in the leader-follower game, the leader (encoder) commits to and announces his optimal policy before the follower (decoder) does, and the follower observes what the leader is committed to before choosing and announcing his optimal policy. Then, a pair of policies $(\mathbf{S}^*, \delta^*)$ is said to be a {\bf Stackelberg equilibrium} \cite{basols99} if
\begin{align}
\begin{split}
&r^t(\mathbf{S}^*, \delta^*(\mathbf{S}^*)) \leq r^t(\mathbf{S}, \delta^*(\mathbf{S})) \quad \forall \,\mathbf{S} \in \mathbb{S}\;, \\
&\hspace{-0.4cm} \text{where } \delta^*(\mathbf{S}) \text{ satisfies} \\
&r^r(\mathbf{S}, \delta^*(\mathbf{S})) \leq r^r(\mathbf{S}, \delta (\mathbf{S})) \quad \forall \,\delta \in \Delta  \,.
\label{eq:stackelbergEquilibrium}
\end{split}
\end{align}
As observed from the definition in \eqref{eq:stackelbergEquilibrium}, the decoder takes his optimal action $\delta^*(\mathbf{S})$ after observing the policy of the encoder $\mathbf{S}$. Further, in the Stackelberg game, the leader cannot backtrack on his commitment, and he has a leadership role since he can manipulate the follower by anticipating follower's actions.

In game theory, Nash (simultaneous game-play) and Stackelberg (sequential game-play) equilibria are drastically different concepts. Both equilibrium concepts find applications depending on the assumptions on the transmitter in view of the commitment conditions. As discussed in \cite{tacWorkArxiv,dynamicGameArxiv}, in the Nash equilibrium case, building on \cite{SignalingGames}, equilibrium properties possess different characteristics as compared to team problems; whereas for the Stackelberg case, the leader agent is restricted to be committed to his announced policy which leads to similarities with team problem setups \cite{CedricWork,akyolITapproachGame,omerHierarchial}. Since there is no such commitment in the Nash setup; the perturbation in the encoder does not lead to a functional perturbation in decoder's policy, unlike the Stackelberg setup. However, in the context of binary signaling, we will see that the distinction is not as sharp as it is in the case of quadratic signaling games \cite{tacWorkArxiv,dynamicGameArxiv}.  

If an equilibrium is achieved when $\mathbf{S}^*$ is non-informative (e.g., $S_0=S_1$) and $\delta^*$ uses only the priors (since the received message is useless), then we call such an equilibrium a {\it non-informative (babbling) equilibrium}.

\subsection{Related Literature}

Standard binary hypothesis testing has been extensively studied over several decades under different setups \cite{Poor,KayDetection}, which can also be viewed as a decentralized control/team problem among an encoder and a decoder who wish to minimize a common cost criterion. However, there exist many scenarios in which the analysis falls within the scope of game theory; either because the goals of the decision makers are misaligned, or because the probabilistic model of the system is not common knowledge among the decision makers. For example, detecting attacks, anomalies, and malicious behavior in network security can be analyzed under the game theoretic perspective \cite{sandberg2015cyberphysical, teixeira2015secure, networkSecurity,cyberSecuritySmartGrid, dan2010stealth}. In this direction, the hypothesis testing and the game theory approaches can be utilized together to investigate attacker-defender type applications \cite{varshneyHypothesisGame,detectionGameAdversary, gupta2010optimal, gupta2012dynamic, basar1985complete}, multimedia source identification problems \cite{sourceIdentification}, and inspection games \cite{Avenhaus1994,inspectionGames,avenhausInspectLeadership}.  

In particular, the binary signaling problem investigated here can be motivated under different application contexts: subjective priors and the presence of a bias in the transmitter cost function when compared with that of the receiver. The former one, decentralized stochastic control with subjective priors, has been studied extensively in the literature \cite{BasTAC85, TeneketzisVaraiya88, CastanonTeneketzis88}. In this setup, players have a common goal but subjective prior information, which necessarily alters the setup from a team problem to a game problem. The latter one is the adaptation of the biased cost function of the transmitter in \cite{SignalingGames} to the binary signaling problem considered here. We discuss these further in the following.


%

\subsection{Two Motivating Setups}

We present two different scenarios that fit into the binary signaling context discussed here and revisit these setups throughout the paper. 

\subsubsection{Subjective Priors}\label{sec:inconPriors}

Suppose that only the beliefs of the transmitter and the receiver about the prior probabilities of hypothesis $\mathcal{H}_0$ and $\mathcal{H}_1$ differ, and the cost values from the perspective of the transmitter and the receiver are the same. Namely, from transmitter's perspective, the priors are $\pi_0^t$ and $\pi_1^t$, whereas the priors are $\pi_0^r$ and $\pi_1^r$ from receiver's perspective, and $C_{ji}=C^t_{ji}=C^r_{ji}$ for $i,j\in\{0,1\}$. 

The setups in decentralized decision making where the priors of the decision makers may be different is an intensely researched area: Among these, \cite{TeneketzisVaraiya88} and \cite{CastanonTeneketzis88} study decentralized decision making with subjective priors, and \cite{BasTAC85} investigates optimal decentralized decision making where the nature of subjective priors converts a team problem into a game problem (see \cite[Section 12.2.3]{YukselBasarBook} for a comprehensive literature review on subjective priors also from a statistical decision making perspective).

\subsubsection{Biased Transmitter Cost Function}\label{sec:biasedCost}


Consider a binary signaling game in which the transmitter encodes a random binary signal $x=i$ as $\mathcal{H}_i$ by choosing the corresponding signal level $S_i$ for $i\in\{0,1\}$, and the receiver decodes the received signal $y$ as $u=\delta(y)$. Let the priors from the perspectives of the transmitter and the receiver be the same; i.e., $\pi_i=\pi_i^t=\pi_i^r$ for $i\in\{0,1\}$, and the Bayes risks of the transmitter and the receiver be defined as $r^t(\mathbf{S},\delta)=\mathbb{E}[\mathds{1}_{\{1=(x\oplus u\oplus b)\}}]$ and $r^r(\mathbf{S},\delta)=\mathbb{E}[\mathds{1}_{\{1=(x\oplus u)\}}]$, respectively, where $\mathds{1}_{\{D\}}$ denotes the indicator function of an event $D$, $\oplus$ stands for the exclusive-or operator and $b$ is the random binary bias term, so that the structure of the costs (Bayes risks) resemble the ones in \cite{SignalingGames} (as also studied in \cite{tacWorkArxiv,dynamicGameArxiv}). Also let $\alpha\triangleq\mathsf{Pr}(b=0)=1-\mathsf{Pr}(b=1)$; i.e., the probability that the cost functions of the transmitter and the receiver are aligned. The following relations can be observed:
\begin{align*}
	r^t(\mathbf{S},\delta)&=\mathbb{E}[\mathds{1}_{\{1=(x\oplus u\oplus b)\}}]\\
	&=\alpha(\pi_0\mathsf{P}_{10}+\pi_1\mathsf{P}_{01})+(1-\alpha)(\pi_0\mathsf{P}_{00}+\pi_1\mathsf{P}_{11}) \\
	& \Rightarrow \quad C^t_{01}=C^t_{10}=\alpha \text{ and } C^t_{00}=C^t_{11}=1-\alpha \;, \\
	r^r(\mathbf{S},\delta)&=\mathbb{E}[\mathds{1}_{\{1=(x\oplus u)\}}]=\pi_0\mathsf{P}_{10}+\pi_1\mathsf{P}_{01} \\
	& \Rightarrow \quad C^r_{01}=C^r_{10}=1 \text{ and } C^r_{00}=C^r_{11}=0 \;.
\end{align*}


\subsection{Contributions}

The main contributions of this study can be summarized as follows:
\begin{enumerate}
	\item A game theoretic formulation of the binary signaling problem is proposed under subjective priors and/or subjective costs. 
	\item Stackelberg and Nash equilibrium policies are obtained and their properties (such as uniqueness and informativeness) are investigated. 
	\begin{enumerate}
		\item It is proved that an equilibrium is almost always informative for a team setup (practically, $0<\tau<1$), whereas in the case of subjective priors and/or costs, it may cease to be informative.
		\item It is shown that Stackelberg equilibria always exist, whereas there are setups under which Nash equilibria may not exist.
	\end{enumerate}
	\item Robustness of equilibrium solutions to small perturbations in the priors or costs are established. It is shown that, the game equilibrium behavior around the team setup is robust under the Nash assumption, whereas it is not robust under the Stackelberg assumption. 
\end{enumerate}

The remainder of the paper is organized as follows. The team setup, the Stackelberg setup, and the Nash setup of the binary signaling game are investigated in Sections II, Section III, and Section IV, respectively. Section V concludes the paper.

\section{TEAM SETUP ANALYSIS}

Now consider the team setup where the cost parameters and the priors are assumed to be same for both the transmitter and the receiver; i.e., $C_{ji}=C^t_{ji}=C^r_{ji}$ and $\pi_i=\pi_i^t=\pi_i^r$ for $i,j\in\{0,1\}$. Thus the common Bayes risk becomes $r^t(\mathbf{S},\delta)=r^r(\mathbf{S},\delta)=\pi_0 (C_{00} \mathsf{P}_{00} + C_{10} \mathsf{P}_{10}) + \pi_1 (C_{01} \mathsf{P}_{01} + C_{11} \mathsf{P}_{11})$. The arguments for the proof of the following result follow from the standard analysis in the detection and estimation literature \cite{Poor,KayDetection}. 
\begin{thm}
	Let $\tau\triangleq{\pi_0 (C_{10}-C_{00}) \over \pi_1 (C_{01}-C_{11})}$. If $\tau\leq0$ or $\tau=\infty$, then the team solution of the binary signaling setup is non-informative. Otherwise; i.e., if $0<\tau<\infty$, the team solution is always informative.
	\label{thm:teamCase}
\end{thm}
\section{STACKELBERG GAME ANALYSIS}

Under the Stackelberg assumption, first the transmitter (the leader agent) announces and commits to a particular policy, and then the receiver (the follower agent) acts accordingly. In this direction, first the transmitter chooses optimal signals $\mathbf{S}=\{S_0,S_1\}$ to minimize his Bayes risk $r^t(\mathbf{S},\delta)$, then the receiver chooses an optimal decision rule $\delta$ accordingly to minimize his Bayes risk $r^r(\mathbf{S},\delta)$. Due to the sequential structure of the Stackelberg game, the encoder knows the priors and the cost parameters of the decoder so that he can adjust his optimal policy accordingly. On the other hand, the decoder knows only the policy of the encoder as he announces during the game-play. Under such a game-play assumption, the equilibrium structure of the Stackelberg binary signaling game can be characterized as follows:

\begin{thm}
If $\tau\triangleq{\pi_0^r (C^r_{10}-C^r_{00}) \over \pi_1^r (C^r_{01}-C^r_{11})}\leq0$ or $\tau=\infty$, the Stackelberg equilibrium of the binary signaling game is non-informative. Otherwise, let $d\triangleq{|S_1-S_0|\over\delta}$, $d_{\max}\triangleq{\sqrt{P_0}+\sqrt{P_1}\over\sigma}$, $\zeta\triangleq\text{sgn}(C^r_{01}-C^r_{11})$, $k_0\triangleq\pi_0^t\zeta(C^t_{10}-C^t_{00})\tau^{-{1\over2}}$, and $k_1\triangleq\pi_1^t\zeta(C^t_{01}-C^t_{11})\tau^{1\over2}$, where the sign of $x$ is defined as 
\[ \text{sgn}(x)=\begin{cases} 
-1 & \text{if }x < 0 \\
0 & \text{if }x=0 \\
1 & \text{if }x>0 
\end{cases}\,.\] Then, the Stackelberg equilibrium structure can be characterized as in Table~\ref{table:stackelbergSummary}, where $d^*=0$ stands for a non-informative equilibrium, and a nonzero $d^*$ corresponds to an informative equilibrium.

\begin{table*}[ht]
	\centering
	\caption{Stackelberg equilibrium analysis for $0<\tau<\infty$.}
	\label{table:stackelbergSummary}
	\begin{adjustbox}{max width=\textwidth}
	\begin{tabular}{|c|c|c|c|}
		\hline
		& \boldmath$\ln\tau\;(k_0-k_1)<0$ & \boldmath$\ln\tau\;(k_0-k_1)\geq0$ \\ \hline
		\boldmath$k_0+k_1<0$	& $d^*=\min\Big\{d_{\max}, \sqrt{\Big|{2\ln\tau(k_0-k_1)\over(k_0+k_1)}\Big|}\Big\}$ & $d^*=0$, non-informative  \\ \hline
		\boldmath$k_0+k_1\geq0$	& $d^*=d_{\max}$ & {$\!\begin{aligned}
			d_{\max}^2<\Big|{2\ln\tau(k_0-k_1)\over(k_0+k_1)}\Big|&\Rightarrow d^*=0\text{, non-informative} \\
			d_{\max}^2\geq\Big|{2\ln\tau(k_0-k_1)\over(k_0+k_1)}\Big|&\Rightarrow \left({k_1\over k_0\tau}\right)^{\text{sgn}(\ln(\tau))}\mathcal{Q}\left({|\ln(\tau)|\over d_{\max}}-{d_{\max}\over2}\right)-\mathcal{Q}\left({|\ln(\tau)|\over d_{\max}}+{d_{\max}\over2}\right)\overset{d^*=d_{\max}}{\underset{d^*=0}{\gtreqless}}0 \end{aligned}$} \\ \hline
	\end{tabular}
	\end{adjustbox}
\end{table*}
\label{thm:stackelbergCases}
\end{thm}

Now we make the following remark on informativeness of the Stackelberg equilibrium:
\begin{remark}
	As we observed in Theorem~\ref{thm:teamCase}, for a team setup, an equilibrium is almost always informative (practically, $0<\tau<\infty$), whereas in the case of subjective priors and/or costs, it may cease to be informative. 
\end{remark}

\begin{figure}[ht]
	\centering
	\includegraphics[width=0.8\linewidth]{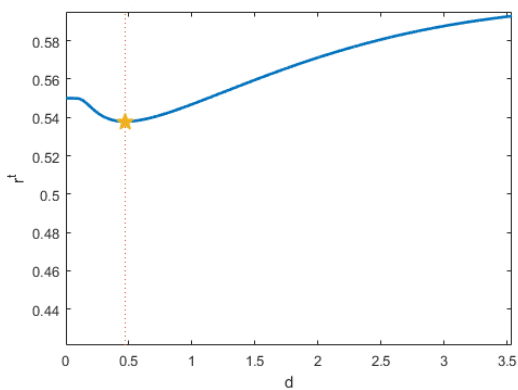}
	\caption{The Bayes risk of the transmitter versus $d$ when $C^r_{01}=0.4, C^r_{10}=0.9, C^r_{00}=0, C^r_{11}=0, C^t_{01}=0.4, C^t_{10}=0.4,  C^t_{00}=0.6, C^t_{11}=0.6, P_0=1,	P_1=1, \sigma = 0.1, \pi_0^t=0.25$, and $\pi_0^r=0.25$. The optimal $d^*=\sqrt{\Big|{2\ln\tau(k_0-k_1)\over(k_0+k_1)}\Big|}=0.4704<d_{max}$ and its corresponding Bayes risk are indicated by the star.}
	\label{fig:StackelbergCase3}
\end{figure}

The most interesting case is when $\ln\tau\;(k_0-k_1)<0, k_0+k_1<0,$ and $d_{\max}^2\geq\Big|{2\ln\tau(k_0-k_1)\over(k_0+k_1)}\Big|$, since in all other cases, the transmitter chooses either the minimum or maximum distance between the signal levels. Further, for classical hypothesis-testing in the team setup, the optimal distance corresponds to maximum separation \cite{Poor}. However, as it can be seen in Figure~\ref{fig:StackelbergCase3}, there is an optimal distance $d^*=\sqrt{\Big|{2\ln\tau(k_0-k_1)\over(k_0+k_1)}\Big|}<d_{\max}$ that makes the Bayes risk of the transmitter minimum.

We now investigate the effects of small perturbations in priors and costs on equilibrium values. In particular, we consider the perturbations around the team setup; i.e., at the point of identical priors and costs. 

Define the perturbation around the team setup as $\vepsilon=\{\epsilon_{\pi0},\epsilon_{\pi1},\epsilon_{00},\epsilon_{01},\epsilon_{10},\epsilon_{11}\}\in\mathbb{R}^6$ such that $\pi_i^t=\pi_i^r+\epsilon_{\pi i}$ and $C_{ji}^t=C_{ji}^r+\epsilon_{ji}$ for $i,j\in\{0,1\}$ (note that the transmitter parameters are perturbed around the receiver parameters which are assumed to be fixed). Then, for $0<\tau<\infty$, at the point of identical priors and costs, small perturbations in both priors and costs imply $k_0=(\pi_0^r+\epsilon_{\pi 0})\zeta(C^r_{10}-C^r_{00}+\epsilon_{10}-\epsilon_{00})\tau^{-{1\over2}}$ and $k_1=(\pi_1^r+\epsilon_{\pi 1})\zeta(C^r_{01}-C^r_{11}+\epsilon_{01}-\epsilon_{11})\tau^{1\over2}$. Since, for $0<\tau<\infty$, $k_0=k_1=\sqrt{\pi^r_0\pi^r_1}\sqrt{(C^r_{10}-C^r_{00})(C^r_{01}-C^r_{11})}>0$ at the point of identical priors and costs, it is possible to obtain both positive and negative $(k_0-k_1)$ by choosing the appropriate perturbation $\vepsilon$ around the team setup. Then, as it can be observed from Table~\ref{table:stackelbergSummary}, even the equilibrium may alter from an informative one to a non-informative one; hence, under the Stackelberg equilibrium, the equilibrium behavior is not robust to small perturbations in both priors and costs.

\subsection{Motivating Examples}
\begin{enumerate}
	\item \textit{Subjective Priors} : Referring to Section~\ref{sec:inconPriors}, the related parameters can be found as follows:
	\begin{align*}
		\tau&={\pi_0^r (C_{10}-C_{00}) \over \pi_1^r (C_{01}-C_{11})} \,,\\
		k_0&=\pi_0^t\sqrt{\pi_1^r\over\pi_0^r}\sqrt{(C_{10}-C_{00})(C_{01}-C_{11})}\,,\\
		k_1&=\pi_1^t\sqrt{\pi_0^r\over\pi_1^r}\sqrt{(C_{10}-C_{00})(C_{01}-C_{11})}\,.
	\end{align*}
	If $0<\tau<\infty$, then $k_0+k_1>0$, and depending on the values of $\ln\tau\;(k_0-k_1)$, $d_{\max}^2$, and $\Big|{2\ln\tau(k_0-k_1)\over(k_0+k_1)}\Big|$, the Stackelberg equilibrium structure can be characterized as in Table~\ref{table:stackelbergInconsistentPriors}.
	Otherwise; i.e., if $\tau\leq0$ or $\tau=\infty$, the equilibrium is non-informative. 
	\begin{table*}[ht]
		\centering
		\caption{Stackelberg equilibrium analysis of subjective priors case for $0<\tau<\infty$.}
		\label{table:stackelbergInconsistentPriors}
		\begin{adjustbox}{max width=\textwidth}
		\begin{tabular}{|c|c|c|c|}
			\hline
			& \boldmath$0<\tau<1$ & \boldmath$1\leq\tau<\infty$ \\ \hline
			\boldmath${\pi_0^t\over\pi_1^t}<{\pi_0^r\over\pi_1^r}$	& {$\!\begin{aligned}
			d_{\max}^2<\Big|{2\ln\tau(k_0-k_1)\over(k_0+k_1)}\Big|&\Rightarrow \text{ Case-5 applies, } d^*=0\text{, non-informative} \\
			d_{\max}^2\geq\Big|{2\ln\tau(k_0-k_1)\over(k_0+k_1)}\Big|&\Rightarrow\text{ Case-6 applies}   \end{aligned}$} & Case-1 applies, $d^*=d_{\max}$  \\ \hline
			\boldmath${\pi_0^t\over\pi_1^t}\geq{\pi_0^r\over\pi_1^r}$	& Case-1 applies, $d^*=d_{\max}$ & {$\!\begin{aligned}
				d_{\max}^2<\Big|{2\ln\tau(k_0-k_1)\over(k_0+k_1)}\Big|&\Rightarrow \text{ Case-5 applies, } d^*=0\text{, non-informative} \\
				d_{\max}^2\geq\Big|{2\ln\tau(k_0-k_1)\over(k_0+k_1)}\Big|&\Rightarrow\text{ Case-6 applies}   \end{aligned}$} \\ \hline
		\end{tabular}
		\end{adjustbox}
	\end{table*}

	\item \textit{Biased Transmitter Cost Function} : Based on the arguments in Section~\ref{sec:biasedCost}, the related parameters can be found as follows:
	\begin{align*}
		\tau={\pi_0\over\pi_1} \,,\; k_0=\sqrt{\pi_0\pi_1}(2\alpha-1)\,,\;
		k_1=\sqrt{\pi_0\pi_1}(2\alpha-1)\,.
	\end{align*}
	Thus, $\ln\tau\;(k_0-k_1)=0$ and $k_0+k_1=2\sqrt{\pi_0\pi_1}(2\alpha-1)$. Then, if $\alpha<1/2$, the transmitter chooses $S_0=S_1$ to minimize $d$, and the equilibrium is non-informative; i.e., he does not send any meaningful information to the transmitter and the receiver considers only the priors. If $\alpha=1/2$, the transmitter has no control on his Bayes risk, hence the  equilibrium is non-informative. Otherwise; i.e., if $\alpha>1/2$, the equilibrium is always informative. In other words, if $\alpha>1/2$, the players act like a team. 
\end{enumerate}

\section{NASH GAME ANALYSIS}

Under the Nash assumption, the transmitter chooses optimal signals $\mathbf{S}=\{S_0,S_1\}$ to minimize $r^t(\mathbf{S},\delta)$, and the receiver chooses optimal decision rule $\delta$ to minimize $r^r(\mathbf{S},\delta)$ simultaneously. In this Nash setup, the encoder and the decoder do not know the priors and the cost parameters of each other; they know only their policies as they announce to each other. Further, there is no commitment between the transmitter and the receiver; hence, the perturbation in the encoder does not lead to a functional perturbation in decoder's policy, unlike the Stackelberg setup. Due to this drastic difference, the equilibrium structure and convergence properties of the Nash equilibrium show significant differences from the ones in the Stackelberg equilibrium, as stated in the following theorem:
\begin{thm}
	Let $\tau\triangleq{\pi_0^r (C^r_{10}-C^r_{00}) \over \pi_1^r (C^r_{01}-C^r_{11})}$ and $\zeta\triangleq\text{sgn}(C^r_{01}-C^r_{11})$, $\xi_0 \triangleq {C^t_{10}-C^t_{00}\over C^r_{10}-C^r_{00}}$, and $\xi_1 \triangleq {C^t_{01}-C^t_{11}\over C^r_{01}-C^r_{11}}$. If $\tau\leq0$ or $\tau=\infty$, then the Nash equilibrium of the binary signaling game is non-informative. Otherwise; i.e., if $0<\tau<\infty$, the Nash equilibrium structure is as depicted in Table~\ref{table:nashSummary}.
	\begin{table*}[ht]
	\centering
	\caption{Nash equilibrium analysis for $0<\tau<\infty$.}
	\label{table:nashSummary}
	\begin{adjustbox}{max width=\textwidth}
		\begin{tabular}{|c|c|c|c|}
			\hline
			& \boldmath$\xi_0>0$ & \boldmath$\xi_0=0$ & \boldmath$\xi_0<0$ \\ \hline
			\boldmath$\xi_1>0$	& unique informative equilibrium & non-informative equilibrium & {$\!\begin{aligned}
				P_0>P_1&\Rightarrow \text{ no equilibrium} \\
				P_0=P_1&\Rightarrow \text{ non-informative equilibrium} \\
				P_0<P_1&\Rightarrow \text{ unique informative equilibrium}
				\end{aligned}$} \\ \hline
			\boldmath$\xi_1=0$	& non-informative equilibrium & non-informative equilibrium & non-informative equilibrium \\ \hline
			\boldmath$\xi_1<0$	& {$\!\begin{aligned}
				P_0>P_1&\Rightarrow \text{ unique informative equilibrium} \\
				P_0=P_1&\Rightarrow \text{ non-informative equilibrium} \\
				P_0<P_1&\Rightarrow \text{ no equilibrium}
				\end{aligned}$} & non-informative equilibrium & no equilibrium \\ \hline
		\end{tabular}
	\end{adjustbox}
\end{table*}
\label{thm:nashCases}
\end{thm}

The main reason for the absence of a non-informative (babbling) equilibrium under the Nash assumption is that in the binary signaling game setup, the receiver is forced to make a decision. Using only the prior information, the receiver always chooses one of the hypothesis. By knowing this, the encoder can manipulate his signaling strategy for his own benefit. However, after this manipulation, the receiver no longer keeps his decision rule the same; namely, the best response of the receiver alters based on the signaling strategy of the transmitter, which entails another change of the best response of the transmitter. Due to such an infinite recursion, there does not exist a pure Nash equilibrium.

Similar to the Stackelberg case, the effects of small perturbations in priors and costs on equilibrium values around the team setup are investigated for the Nash setup as follows:

Define the perturbation around the team setup as $\vepsilon=\{\epsilon_{\pi0},\epsilon_{\pi1},\epsilon_{00},\epsilon_{01},\epsilon_{10},\epsilon_{11}\}\in\mathbb{R}^6$ such that $\pi_i^t=\pi_i^r+\epsilon_{\pi i}$ and $C_{ji}^t=C_{ji}^r+\epsilon_{ji}$ for $i,j\in\{0,1\}$  (note that the transmitter parameters are perturbed around the receiver parameters which are assumed to be fixed). Then, for $0<\tau<\infty$, at the point of identical priors and costs, small perturbations in priors and costs imply $\xi_0 = {C^r_{10}-C^r_{00}+\epsilon_{10}-\epsilon_{00}\over C^r_{10}-C^r_{00}}$ and $\xi_1 = {C^r_{01}-C^r_{11}+\epsilon_{01}-\epsilon_{11}\over C^r_{01}-C^r_{11}}$. As it can be seen, the Nash equilibrium is not affected by small perturbations in priors. Further, since $\xi_0=\xi_1=1$ at the point of identical priors and costs for $0<\tau<\infty$, as long as the perturbation $\vepsilon$ is chosen such that $\lvert{\epsilon_{10}-\epsilon_{00}\over C^r_{10}-C^r_{00}}\rvert<1$ and $\lvert{\epsilon_{01}-\epsilon_{11}\over C^r_{01}-C^r_{11}}\rvert<1$, we always obtain positive $\xi_0$ and $\xi_1$ in Table~\ref{table:nashSummary}. Thus, under the Nash assumption, the equilibrium behavior is robust to small perturbations in both priors and costs.

\subsection{Motivating Examples}
\begin{enumerate}
	\item \textit{Subjective Priors} : The related parameters are $\tau={\pi_0^r (C_{10}-C_{00}) \over \pi_1^r (C_{01}-C_{11})}$, $\xi_0 =1$, and $\xi_1=1$. Thus, if $\tau<0$ or $\tau=\infty$, the equilibrium is non-informative; otherwise, there always exists a unique informative equilibrium; namely, as long as the priors are mutually absolutely continuous, the subjectivity in the priors does not affect the equilibrium. 
	\item \textit{Biased Transmitter Cost Function} : If $\alpha>1/2$, the equilibrium is informative; if $\alpha=1/2$, the equilibrium is non-informative; otherwise; i.e., if $\alpha<1/2$, there exists no equilibrium. As it can be seen, the existence of the equilibrium depends on $\alpha=\mathsf{Pr}(b=0)$, the probability that the Bayes risks of the transmitter and the receiver are aligned. 
\end{enumerate}

\section{CONCLUSION}

In this paper, we considered binary signaling problems in which the decision makers (the transmitter and the receiver) have subjective priors and/or misaligned objective functions. Depending on the commitment nature of the transmitter to his policies, we formulated the binary signaling problem as a Bayesian game under either Nash or Stackelberg equilibrium concepts and established equilibrium solutions and their properties. We showed that there can be informative or non-informative equilibria in the binary signaling game under the Stackelberg assumption, but there always exists an equilibrium. However, apart from the informative and non-informative equilibria cases, there may not be a Nash equilibrium when the receiver is restricted to use deterministic policies. We also studied the effects of small perturbations at the point of identical priors and costs and showed that the game equilibrium behavior around the team setup is robust under the Nash assumption, whereas it is not robust under the Stackelberg assumption.

%


\bibliographystyle{IEEEtran}
\bibliography{../SerkanBibliography}

\end{document}